%% file: main.tex
\documentclass[12pt, reqno]{amsart}
\usepackage{amssymb,latexsym, mathtools,tikz, hyperref}
\usetikzlibrary{arrows,decorations.markings}
\tikzstyle arrowstyle=[scale=1]
\tikzstyle directed=[postaction={decorate,
decoration={markings,mark=at position .65 with {\arrow[arrowstyle]{stealth}}}}]
\usepackage{xypic}
\usepackage[top=35mm, bottom=25mm, left=30mm, right=30mm]{geometry}

\usepackage{microtype}
\usepackage{enumerate}
\usepackage{graphicx}
\usepackage{float}
\usepackage{placeins}
\usepackage{mdframed}
\usepackage{amssymb}
\usepackage{esint}
\usepackage{cool}
\usepackage[all,cmtip]{xy}
\usepackage{mathtools}
\usepackage{amstext} 
\usepackage{array}   
\usepackage[shortlabels]{enumitem}
\usepackage{ytableau}
\usepackage[ruled,vlined]{algorithm2e}

\usepackage{cleveref}
\usepackage[utf8]{inputenc}
\usepackage{comment}
\input macro.tex

\begin{document}

\title[Linear Strands of Initial Ideals of DFIs]{Linear Strands of Initial Ideals of Determinantal Facet Ideals}

\author{Ayah Almousa}
\address{University of Minnesota - Twin Cities}
\email{almou007@umn.edu}
\urladdr{\url{http://umn.edu/~almou007}}

\author{Keller VandeBogert}
\address{University of Notre Dame}
\email{kvandebo@nd.edu}
\urladdr{\url{https://sites.google.com/view/kellervandebogert/}}

\keywords{determinantal facet ideal, binomial edge ideal, initial ideals, linear strand, free resolutions, cellular resolutions}
\subjclass[2020]{Primary: 13F20,13F55; Secondary: 55U10,05E40}
\date{\today}
\begin{abstract}
We construct an explicit linear strand for the initial ideal of any determinantal facet ideal (DFI) with respect to any diagonal term order. We show that if the clique complex of $\Delta$ has no \emph{1-nonfaces} larger than a certain cardinality, then the Betti numbers of the linear strand of $J_\Delta$ and its initial ideal coincide. This confirms a conjecture of Ene--Herzog--Hibi for closed graphs with at most $2$ maximal cliques. Additionally, we show that the linear strand of the initial ideal of any DFI is supported on an induced subcomplex of the complex of boxes introduced by Nagel--Reiner.
\end{abstract}
\maketitle

\input sections/intro

\input sections/Background

\input sections/SparseEN
\input sections/cellular

\section*{Acknowledgments}
The authors would like to thank the anonymous referee for their close reading and helpful comments. The first author was partially supported by the NSF GRFP under Grant No. DGE-1650441.

\bibliographystyle{amsplain}
\bibliography{biblio}
\addcontentsline{toc}{section}{Bibliography}

\end{document}

%% file: macro.tex
\newcolumntype{L}{>{$}l<{$}} 
\newtheorem{lemma}{Lemma}[section]
\newtheorem{theorem}[lemma]{Theorem}

\newtheorem{prop}[lemma]{Proposition}
\newtheorem{cor}[lemma]{Corollary}
\newtheorem{conj}[lemma]{Conjecture}

\theoremstyle{remark}
\newtheorem{remark}[lemma]{Remark}

\theoremstyle{definition}
\newtheorem{definition}[lemma]{Definition}
\newtheorem{example}[lemma]{Example}
\newtheorem{setup}[lemma]{Setup}
\newtheorem{notation}[lemma]{Notation} 
\newtheorem{construction}[lemma]{Construction}

\newcommand{\lin}{\operatorname{lin}}

\newcommand{\reg}{\operatorname{reg}}

\newcommand{\inn}{\operatorname{in}}
\newcommand{\w}{\wedge}

\renewcommand{\I}{\mathcal{I}}

\newcommand{\cat}[1]{\mathcal{#1}}
\newcommand{\clique}{\Delta^{\textrm{clique}}}
\newcommand{\del}{\partial}
\newcommand{\one}{\mathbbm{1}}

\newcommand{\coker}{\operatorname{coker}}

\renewcommand{\hom}{\operatorname{Hom}} 

\newcommand{\pd}{\operatorname{pd}}

\newcommand{\mdeg}{\operatorname{mdeg}}

\newcommand{\doot}{\bullet}

\newcommand{\ff}{\mathbf f}

\renewcommand{\aa}{\mathbf a}
\newcommand{\bb}{\mathbf b}
\newcommand{\cc}{\mathbf c}

\newcommand{\xx}{\mathbf x}

\newcommand{\cF}{\mathcal{F}}

\newcommand{\cP}{\mathcal{P}}

\newcommand{\cX}{\mathcal{X}}

\newcommand{\NN}{\mathbb{N}}

\newcommand{\ZZ}{\mathbb{Z}}

\newcommand{\Xv}{\check{X}}

\newcommand{\ra}{\rightarrow}
\renewcommand{\del}{\partial}

\DeclarePairedDelimiter\abs{\lvert}{\rvert}%
\DeclareMathOperator{\lcm}{lcm}

%% file: sections/intro.tex
\section{Introduction}

Let $R$ denote the coordinate ring of a generic $n \times m$ matrix $M$, over some field $k$ (with $n \leq m$). The ideal of maximal minors $I_n (M)$ possesses many surprising and desirable properties; for example, a result of Sturmfels and Zelevinsky \cite{SturmfelsZelevinsky93} shows that the natural generating set consisting of all maximal minors of $M$ forms a reduced Gr\"obner basis of $I_n (M)$ for \emph{any} term order $<$. Boocher goes one step further in \cite{boocher2012} and shows that the graded Betti numbers of $I_n (M)$ and any of its initial ideals must also agree. As a consequence, the minimal free resolution of any initial ideal of the ideal of maximal minors can be obtained by simply setting some of the entries in the matrix representation of the standard Eagon-Northcott differentials equal to $0$. 

One direction for generalizing ideals of maximal minors is to imagine that the column-sets of minors appearing in the generating set are parametrized by some simplicial complex $\Delta$. Such ideals are called \emph{determinantal facet ideals} (DFIs), and have been studied in multiple contexts by a wide variety of authors (see \cite{ene2011determinantal}, \cite{herzog2015linear}, \cite{vdb2020}). In the case that $\Delta$ has dimension $1$ (that is, $\Delta$ is a graph), DFIs are better known as \emph{binomial edge ideals}, and are even more well behaved than arbitrary DFIs (see \cite{ohtani2011}, \cite{herzog2010binomial}, and \cite{madani16-BEIsurvey}). DFIs themselves have also been generalized in a few different directions - Mohammadi and Rauh \cite{mohammadi2018prime} have allowed for the minors to be parametrized by hypergraphs, and in \cite{almousa2020resolutions}, DFIs for nonmaximal minors are studied. 

In \cite{EHH2011-CMBEIs}, Ene, Herzog, and Hibi conjecture that Betti numbers for a \emph{closed} binomial edge ideal agree with that of its initial ideal with respect to any diagonal term order. This conjecture is known to be true in the case of Cohen-Macaulay binomial edge ideals, but has remained elusive in generality. In this paper, we give further evidence and a proof of this conjecture in the $2$-clique case using techniques related to the computation of linear strands.

Given a homogeneous minimal complex of initial degree $d$, the linear strand can be obtained by restricting the $i$th differential to basis elements of degree $i+d$, for all $i \geq 1$. The first results relating to linear strands of DFIs were given by Herzog, Kiani, and Madani in \cite{herzog2015linear}; the main result was the fact that the linear strand of any DFI can be obtained as a so-called \emph{generalized} Eagon-Northcott complex. In this paper, we prove a similar result showing that the linear strand of the initial ideal of certain classes of DFIs is obtained as a \emph{generalized sparse} Eagon-Northcott complex (see Definition \ref{def:gensparseEN}). 

More precisely, in this paper we study linear strands of initial ideals of DFIs with respect to a diagonal term order. We first construct an explicit resolution of $\inn_< I_n (M)$ as a \emph{sparse} Eagon-Northcott complex (see \ref{def:resnoflexideal}). Using this, one can restrict to an appropriate subcomplex parametrized by the clique complex of the simplicial complex associated to the DFI $J_\Delta$. This subcomplex is not acyclic in general; however, if one imposes sufficient conditions on $\Delta$, the homology of this subcomplex will vanish in ``small" degrees. This implies that the Betti numbers of certain DFIs and their initial ideals with respect to a diagonal term order coincide on the linear strand (see Theorem \ref{thm:mainpunchline}). 

As an application, we prove the previously mentioned conjecture of Ene, Herzog, and Hibi in the case that the associated graph $G$ has at most $2$ maximal cliques. Moreover, we pose more generally the conjecture that for any lcm-closed DFI $J_\Delta$, the Betti numbers of $J_\Delta$ and its initial ideal with respect to any diagonal term order coincide (see Conjecture \ref{conj:BnosConj}). Finally, we consider the linear strands of initial ideals (with respect to any diagonal term order) of arbitrary DFIs. In this case, we find that the linear strand in general is \emph{always} supported on a polyhedral cell complex. This implies that the Betti numbers of a general linear strand can be obtained by looking at certain induced subcomplexes of the so-called \emph{complex of boxes} of Nagel and Reiner (see \cite{NR09}). We conclude with examples of this construction and remarks on further applications.

The paper is organized as follows. In Section \ref{sec:Background}, we introduce determinantal facet ideals and notation and terminology surrounding them that will be in play for the rest of the paper.
We open Section \ref{sec:sparseEN} with the definition of the Eagon-Northcott complex and the aforementioned result of Boocher on the Betti numbers of initial ideals of ideals of maximal minors.
We go on to introduce a \emph{sparse} Eagon-Northcott complex to be used for resolving the initial ideal of the ideal of maximal minors with respect to any diagonal term order. As corollaries, we obtain minimal free resolutions of the ideal of squarefree monomials of a given degree and the box polarization for powers of the graded maximal ideal (and hence by specialization, any power of the graded maximal ideal).

In Section \ref{sec:linStrandlcmclosed}, we consider subcomplexes of the sparse Eagon-Northcott complex with basis elements parametrized by a simplicial complex $\Delta$. This subcomplex in general is \emph{not} acyclic, but it turns out that combinatorial properties of $\Delta$ will allow us to deduce exactly when homology is nontrivial in linear degrees. This combinatorial condition is encoded in the existence of so-called $i$-\emph{nonfaces}, a generalization of minimal nonfaces (see Definition \ref{def:inonface}). We use this to show that the Betti numbers along the linear strand of a DFI $J_\Delta$ and its initial ideal with respect to a diagonal term order coincide if and only if the clique complex associated to $\Delta$ has no $1$-nonfaces of cardinality $n+1$. We apply this result to give a proof of the conjecture of Ene, Herzog, and Hibi for graphs having at most $2$ maximal cliques (see Corollary \ref{cor:EHHconj}).

In Section \ref{sec: cellular}, we consider linear strands supported on cellular complexes. In particular, we show that the linear strand of the initial ideal (with respect to a diagonal term order) of \emph{any} DFI is supported on a polyhedral cellular complex which can be obtained as the induced subcomplex of the \emph{complex of boxes} introduced by Nagel and Reiner in \cite{NR09} (see Theorem \ref{thm: cplxBoxesLinStrand}). This implies that the multigraded Betti numbers of the linear strand of any such ideal are $0$ or $1$ and moreover allows us to count the Betti numbers in the linear strand as the $f$-vector of the associated polyhedral cell complex.

%% file: sections/Background.tex
\section{Background}\label{sec:Background}

In this section we introduce some necessary background to be used for the rest of the paper. In particular, we discuss determinantal facet ideals and establish some notation related to these ideals and the simplicial complexes that parametrize them. 

\begin{notation}
Let $R = k[x_{ij} \mid 1 \leq i \leq n, \ 1 \leq j \leq m]$ be a polynomial ring over an arbitrary field $k$ and let $M = (x_{ij})_{1\leq i \leq n, 1 \leq j \leq m}$ denote a generic $n \times m$ matrix where $n \leq m$.
For indices $\aa = \{a_1,\ldots, a_r\}$ and $\bb = \{b_1, \ldots, b_r\}$ such that $1\leq a_1 < \ldots < a_r\leq n$ and $1\leq b_1 < \cdots < b_r\leq m$, set
$$
[\aa\vert\bb] = [a_1,\ldots, a_r\vert b_1,\ldots, b_r] = \det \left( \begin{array}{ccc}
    x_{a_1,b_1} & \cdots & x_{a_1,b_r} \\
    \vdots & \ddots & \vdots\\
    x_{a_r,b_1} & \cdots & x_{a_r,b_r}\\
\end{array} \right)
$$
where $[\aa\vert\bb]=0$ if $r>n$. When $r=n$, use the simplified notation $[\aa]$ = $[1,\ldots, n\vert \aa]$. The ideal generated by the $r$-minors of $M$ is denoted $I_r(M)$.
\end{notation}

\begin{definition}\label{def:cliques}
Let $\Delta$ be an $(n-1)$-dimensional simplicial complex on vertex set $[m]$.
A \textit{clique} of $\Delta$ is an induced subcomplex $\Gamma$ of $\Delta$ such that any $n$ vertices of $\Gamma$ are in a face together.
A clique is called \textit{maximal} if it is not contained in any larger clique of $\Delta$.
The simplicial complex $\clique$ whose facets are the maximal cliques of $\Delta$ is called the \emph{clique complex} associated to $\Delta$. The decomposition $\Delta = \Delta_1 \cup \cdots \cup \Delta_c$ is called the \emph{maximal clique decomposition} of $\Delta$.
\end{definition}

\begin{definition}\label{def: DFI}
Let $\Delta$ be a pure $(n-1)$-dimensional simplicial complex on the vertex set $[m]$. Let $R = k[x_{ij} \mid 1 \leq i \leq n, \ 1 \leq j \leq m]$ be a polynomial ring over an arbitrary field $k$ and let $M = (x_{ij})_{1\leq i \leq n, 1 \leq j \leq m}$ denote a generic $n \times m$ matrix where $n \leq m$. The \emph{determinantal facet ideal} (or \emph{DFI}) $J_\Delta\subseteq S$ associated to $\Delta$ is the ideal generated by $n$-minors $[\aa ]$ where $\aa$ supports an $(n-1)$-facet of $\Delta$; that is, the columns of $[\aa ]$ correspond to the vertices of some facet of $ \Delta$.
\end{definition}

\begin{example}
If $\Delta$ is an $(n-1)$-simplex, then $J_\Delta$ is generated by a single polynomial corresponding to the determinant of an $n\times n$ generic matrix.
More generally, when $\Delta$ is the $n$-skeleton of an $(m-1)$-simplex, $J_\Delta$ is the ideal of maximal minors $I_n(M)$.
\end{example}

\begin{notation}\label{not: leadTerm} Let $<$ be a total monomial order in a polynomial ring $R$ over a field $k$. If $f$ is a polynomial in $R$, denote by $\inn_<(f)$ the leading term of $f$ with respect to $<$. If $I\subseteq R$ is an ideal, denote by $\inn_<(I)$ the initial ideal of $I$ with respect to $<$. Frequently, when the term order $<$ is clear, it will be dropped and leading terms and initial ideals will simply be denoted by $\inn(f)$ and $\inn(I)$, respectively.
\end{notation}

Recall that a term order $<$ on $R$ is a \textit{diagonal term order} if the leading term with respect to $<$ of the determinant of a submatrix of $M$ is the product of the indeterminates on its main diagonal.

\begin{notation}
Let $\Delta$ be a pure $(n-1)$-dimensional simplicial complex on the vertex set $[m]$ with maximal clique decomposition $\Delta = \Delta_1 \cup \cdots \cup \Delta_c$. The notation $J_{\Delta_i}$ will be used to denote the DFI associated to the simplicial complex whose facets come from all $\aa \in \Delta_i$ with $|\aa| = n$. 
\end{notation}

\begin{definition}\label{def: lcmClosed} Let $\Delta$ be a pure $(n-1)$-dimensional simplicial complex on $m$ vertices with maximal clique decomposition $\Delta = \bigcup_{i=1}^c \Delta_i$. Let $<$ denote any total monomial order for the polynomial ring $R$. The DFI $J_\Delta$ is \emph{lcm-closed} if the following condition holds:
\begin{enumerate}[(*)]
    \item For all $[\aa ] \in J_{\Delta_i}$ and $[\aa' ] \in J_{\Delta_j}$ with non-coprime lead terms such that $[\aa ], \ [\aa' ] \notin J_{\Delta_i \cap \Delta_j}$, there exists some maximal minor $[\cc ] \in J_{\Delta_i \cap \Delta_j}$ such that $\inn ([\cc ])$ divides $\lcm\big( \inn([\aa ] ) , \inn ( [\aa' ]) \big)$. 
\end{enumerate}
\end{definition}

\begin{example} Let $\Delta$ be a $2$-dimensional simplicial complex on $5$ vertices with clique decomposition $\Delta_1 = \{1,2,3,4\}$ and $\Delta_2 = \{2,3,4,5\}$, and so $n = 3$ and $m = 5$. Let $<$ be any diagonal term order. For any $2$-dimensional faces $\aa\in\Delta_1$ and $\aa' \in \Delta_2$ such that $\aa,\aa'\notin \Delta_1\cap\Delta_2$, observe that $\inn_<[\aa]$ and $\inn_<[\aa']$ are coprime except for the case where $\aa = [1,3,4]$ and $\aa' = [2,3,5]$. In this case,
\begin{align*}
\lcm(\inn_<([1,3,4]), \inn_<([2,3,5]) &= \lcm(x_{11} x_{23} x_{34}, x_{12} x_{23} x_{35}) \\
&= x_{11} x_{12} x_{23} x_{34} x_{35}
\end{align*}
which is divisible by $x_{12} x_{23} x_{34} = \inn_<([2,3,4])$, and $\inn_<([2,3,4]) = J_{\Delta_1}\cap J_{\Delta_2}$. Therefore, $J_\Delta$ is lcm-closed.
\end{example}

\begin{example}
Let $\Delta$ be the $2$-dimensional simplicial complex on $4$ vertices with facets $\Delta_1 = \{1,3,4\}$ and $\Delta_2 = \{2,3,4\}$, so $n =3$ and $m = 4$. Observe that $\Delta_1\cap\Delta_2$ is too small to contain any generators of $\Delta$. Let $<$ be any diagonal term order. Then
\begin{align*}
\lcm(\inn_<([1,3,4]),\inn_<([2,3,4])) = x_{2,3} x_{3,4}
\end{align*}
which is not divisible by any generator of $\Delta$. Therefore, $J_\Delta$ is not lcm-closed.
\end{example}

In \cite{almousa2020resolutions}, it is shown that the standard minimal generating set of an lcm-closed DFI forms a reduced Gr\"obner basis; conjecturally, we believe that Definition \ref{def: lcmClosed} is equivalent to being a Gr\"obner basis for DFIs.

%% file: sections/SparseEN.tex
\section{Sparse Eagon-Northcott Complexes}\label{sec:sparseEN}

We open this section with some results on ideals generated by maximal minors, including the definition of the classical \emph{Eagon-Northcott} complex. This includes a key observation stating that the initial ideal of the ideal of maximal minors with respect to any term order has multigraded Betti numbers that are either $0$ or $1$.
At the heart of this section is the construction of an explicit example of a \emph{sparse} Eagon-Northcott complex. The most complicated part of the construction ends up being the definition of the differentials; as it turns out, ayclicity will follow immediately from Theorem \ref{thm:boocherBnos}. As a consequence, we deduce that certain specializations of this complex also yield minimal free resolutions of the ideal generated by all squarefree monomials of a given degree and powers of the graded maximal ideal in a polynomial ring. We begin this section with the following setup:

\begin{setup}\label{set:linearStrandSetup}
Let $R= k[x_{ij} \mid 1\leq i \leq n, 1 \leq j \leq m]$ and $M = (x_{ij})_{1\leq i \leq n, 1 \leq j \leq m}$ denote a generic $n \times m$ matrix, where $n \leq m$. View $M$ as a homomorphism $M : F \to G$ of free modules $F$ and $G$ of rank $m$ and $n$, respectively, where $\{f_i\}_{i\in [m]}$ and $\{g_j\}_{j\in [n]}$ denote the standard bases with respect to which $M$ has the above matrix representation. Let $<$ denote any diagonal term order on $R$, and denote by $\inn_< I_n(M)$ the initial ideal with respect to $<$ of the ideal of maximal minors of $M$.
\end{setup}

\begin{notation}
Let $R$ be a polynomial ring over a field $k$. Let $F$ be a free $R$-module of rank $m$ with basis $f_1,\dots, f_m$. Denote by $D_d(F)$ the $d$th divided power of $F$, and by $\bigwedge^d F$ the $d$th exterior power of $F$. Let $J = \{ j_1 < j_2 < \dots < j_k\} \subset [m]$. Define
$$
f_J \coloneqq f_{j_1} \wedge \dots \wedge f_{j_k}\in \bigwedge^k F.
$$
If $\alpha = (\alpha_1, \dots, \alpha_n) \in \NN^n$ such that $\abs{\alpha} = d$, set
$$
f^{(\alpha)} \coloneqq f_1^{(\alpha_1)} f_2^{(\alpha_2)} \dots f_n^{(\alpha_n)} \in D_d(F).
$$
\end{notation}

\begin{definition}[Eagon-Northcott complex]\label{def:EN}
Let $\phi : F \to G$ be a homomorphism of free modules of ranks $m$ and $n$, respectively, with $n \leq m$. Let $c_\phi$ be the image of $\phi$ under the isomorphism $\hom_R (F,G) \xrightarrow{\cong} F^* \otimes G $. 
The \emph{Eagon-Northcott complex} is the complex
$$0 \to D_{n-m} (G^*) \otimes \bigwedge^n F \to D_{n-m-1} (G^*) \otimes \bigwedge^{n-1} F \to \cdots \to G^* \otimes \bigwedge^{m+1} F \to \bigwedge^m F \to \bigwedge^m G$$
with differentials in homological degree $\geq 2$ induced by multiplication by the element $c_\phi \in F^* \otimes G$, and the map $\bigwedge^g F \to \bigwedge^g G$ is $\bigwedge^g \phi$. 
\end{definition}

\begin{theorem}[{\cite[Theorem A2.10]{eisenbud2013commutative}}] The Eagon-Northcott complex of Definition \ref{def:EN} is a homogeneous minimal free resolution of $\coker(\bigwedge^m \phi)$.
\end{theorem}

\begin{notation}\label{not:ENgrading}
Let $E_\bullet$ denote the Eagon-Northcott complex of Definition \ref{def:EN}. If $F$ has basis $f_1 , \dots , f_m$ and $G$ has basis $g_1 , \dots , g_n$, then define
$$g^{*(\alpha)} \otimes f_I := g_1^{*(\alpha_1)} \cdots g_n^{*(\alpha_n)} \otimes f_{i_1} \w \cdots \w f_{i_{n+\ell}},$$
where $\alpha = (\alpha_1 , \dots , \alpha_n)\in \NN^n$ and $I = (i_1 < \cdots < i_{n+\ell})$ where $0\leq \ell \leq m-n$. Observe that $E_\bullet$ is $\ZZ^n \times \ZZ^m$-graded via
\begin{equation*}
\mdeg (g^{* (\alpha)} \otimes f_I ) = (\one + \alpha_1 \epsilon_1 + \cdots + \alpha_n \epsilon_n , \epsilon_{i_1} + \cdots +\epsilon_{i_{n+\ell}}),
\end{equation*}
where $\epsilon_k$ denotes the appropriately sized vector with $1$ in the $i$th spot and $0$ elsewhere, and $\one$ denotes a length $n$ vector of $1$s.
\end{notation}

The following result, due to Boocher, shows that with respect to any term order $<$, the ideal $\inn_< I_n (M)$ specializes to the ideal of all squarefree monomials of degree $n$ in $m$ variables. 

\begin{theorem}[{\cite[Proof of Theorem 3.1]{boocher2012}}]\label{thm:boocherBnos}
For any term order $<$, the sequence of variable differences
    $$\{ x_{11} - x_{21} , \dots, x_{11} - x_{n1} \} \cup \cdots \cup \{ x_{1m} - x_{2m} , \dotsc , x_{1m} - x_{nm} \}$$
    forms a regular sequence on $R/ \inn_< I_n (M)$. In particular, all the graded Betti numbers of $I_n(M)$ and $\inn_< I_n(M)$ coincide, i.e.
    $$\beta_{ij} (R / I_n(M)) = \beta_{ij} (R/ \inn_< I_n (M)) \ \textrm{for all} \ i,j.$$
\end{theorem}

\begin{cor}
For every multidegree $\alpha\in \ZZ^{nm}$ and any term order $<$,
$$\beta_{\alpha} (R / \inn_< I_n(M)) \leq 1.$$
\end{cor}

\begin{proof}
By Theorem \ref{thm:boocherBnos}, a minimal free resolution of $R/\inn_< I_n (M)$ may be obtained by setting some of the entries in the matrix representation of the differentials of Definition \ref{def:EN} equal to $0$.
If $\alpha$ is a multidegree with respect to the $\ZZ^n \times \ZZ^m$-grading of the Eagon-Northcott complex $E_\bullet$, one has 
$$\beta_\alpha (R / \inn_< I_n (M)) \leq 1.$$
Let $F_\doot$ be a $\ZZ^{nm}$-graded minimal free resolution of the monomial ideal $\inn_< I_n(M)$.
Then $F_\doot$ is also $\ZZ^n\times \ZZ^m$-graded, where every multidegree $\alpha'\in \ZZ^{nm}$ corresponds to a unique basis element in $E_\doot$ with multidegree $\alpha\in \ZZ^n\times \ZZ^m$ via
$$
x_{1,i_{11}} x_{1,i_{12}} \dots x_{1,i_{1k_1}} x_{2,i_{21}} \dots x_{2,i_{2k_2}} \dots x_{n, i_{n k_n}} \mapsto g^{*(k_1, k_2,\dots, k_n)} \otimes f_I
$$
where $I$ lists the elements of the set $\{i_{11}, i_{12}, \dots, i_{1k_1}, i_{21}, \dots, i_{n k_n}\}$ in increasing order. Therefore $\beta_{\alpha'} (R / \inn_< I_n (M))= \beta_\alpha (R / \inn_< I_n (M))$ and the result follows.
\end{proof}



\begin{notation}
Let $\alpha = (\alpha_1 , \dotsc , \alpha_n)\in \NN^n$. Define
$$\alpha_{\leq i } := (\alpha_1 , \dotsc , \alpha_i ),$$
where $\alpha_{\leq i} = \varnothing$ if $i \leq 0$ and $\alpha_{\leq i } = \alpha$ if $i \geq n$.
\end{notation}

\begin{definition}\label{def:sparseENdiff}
Let $\alpha = (\alpha_1 , \dotsc , \alpha_n)\in \NN^n$ with $|\alpha | = \ell$ and $I = (i_1 < \cdots < i_{n+\ell})$ where $0\leq \ell \leq m-n$. Define the indexing set
$$\I_< ( \alpha , I) := \{ (i , I_{i+j}) \mid i \in \{ k \mid \alpha_k >0 \}, \ |\alpha_{\leq i-1}| \leq j \leq |\alpha_{\leq i}| \}$$
\end{definition}

\begin{example}
One easily computes:
$$\I_< ( (1,1,1) , (1,2,3,4,5,6) ) = \{(1,1), (1,2), (2,3), (2,4), (3,5), (3,6)\}$$
$$\I_< ( (1,0,2),(1,2,3,4,5,6)) = \{(1, 1), (1, 2), (3, 4), (3, 5), (3, 6)\}$$
$$\I_< ((2,1),(1,2,4,5,6)) = \{(1, 1), (1, 2), (1, 4), (2, 5), (2, 6)\}$$
\end{example}

\begin{remark}
Each basis element of the Eagon-Northcott complex also has a $\ZZ^{nm}$-grading and can be viewed as a monomial in the ideal $\inn(J_\Delta)$:
$$
g^{* ( \alpha)} \otimes f_\sigma \leftrightarrow (x_{1\sigma_1} \dots x_{1\sigma_{\alpha_1+1}})\cdot (x_{2\sigma_{\alpha_1+2}}\dots x_{2\sigma_{\alpha_1+\alpha_2+2}}) \cdots (x_{n\sigma_{n+i-\alpha_n-2}}\dots x_{n\sigma_{n+i-1}}) =: m_{\alpha , \sigma}.
$$
Observe then that $\I_< (\alpha , \sigma)$ chooses precisely those indices for which $m_{\alpha , \sigma} / x_{r s} \in \inn_< I_n (M)$ for all $(r,s) \in \I_< (\alpha , \sigma)$.
\end{remark}

Using the indexing set of Definition \ref{def:sparseENdiff}, we now define the sparse Eagon-Northcott complex.

\begin{definition}\label{def:resnoflexideal}
Adopt notation and hypotheses as in Setup \ref{set:linearStrandSetup}. Let $E_\bullet'$ denote the sequence of module homomorphisms with
$$E_\ell' = \begin{cases} 
\bigwedge^n G \qquad \textrm{if} \ \ell = 0 \\
D_\ell (G^*) \otimes \bigwedge^{n+\ell} F \qquad \textrm{otherwise}, \\
\end{cases}$$
and first differential $d'_1 : \bigwedge^n F \to \bigwedge^n G$ sending $f_I \mapsto \inn_<  (M(f_I))$, that is, $f_I$ maps to the leading term of the image of $f_I$ under the homomorphism $M$. For $\ell\geq 2$, $d'_\ell : D_{\ell - 1} (G^*) \otimes \bigwedge^{n+\ell-1} F \to D_{\ell-2} (G^*) \otimes \bigwedge^{n+\ell-2} F$ is the map
$$d_\ell ( g^{*(\alpha)} \otimes f_I ) = \sum_{\{i \mid \alpha_i >0\}} \sum_{(i,I_j) \in \I_< (\alpha , I)} (-1)^{j + 1}  x_{iI_j} g^{*(\alpha - \epsilon_i)} \otimes f_{I \backslash I_j }.$$
\end{definition}

\begin{prop}\label{prop:isacomplex}
The sequence of homomorphisms $E'_\bullet$ of Definition \ref{def:resnoflexideal} forms a complex.
\end{prop}

\begin{proof}
Observe first that the map $d'_1: \bigwedge^n F \to \bigwedge^n G$ sends
$$f_I \mapsto x_{1I_1} \cdots x_{n I_n} g_{[n]}.$$
We first verify that $d'_1 \circ d'_2 = 0$. Let $g_k^* \otimes f_I \in G^* \otimes \bigwedge^{n+1} F$; then:
\begingroup\allowdisplaybreaks
\begin{align*}
    d'_1 \circ d'_2 (g_k^* \otimes f_I ) &= d'_1 ((-1)^{k+1} x_{k I_k} f_{I \backslash I_k} + (-1)^{k+2}x_{k I_{k+1}} f_{I \backslash I_{k+1}} ) \\
    &= (-1)^{k+1} x_{k I_k} (x_{1I_1} \cdots \widehat{x_{k I_k}} x_{k I_{k+1}} \cdots x_{nI_n} )g_{[n]} \\
    &\quad + (-1)^{k+2} x_{k I_{k+1}} (x_{1I_1} \cdots x_{k I_k} \widehat{x_{k I_{k+1}}} \cdots x_{nI_n} )g_{[n]} \\
    &=0. \\
\end{align*}
\endgroup
Assume now that $\ell \geq 1$; the fact that $d'_{\ell+1} \circ d'_{\ell+2} = 0$ is a nearly identical computation to that of the standard Eagon-Northcott differential, where one uses the fact that
$$\I_< ( \alpha - \epsilon_i , I \backslash I_j ) = \begin{cases}
\I_< ( \alpha , I) \backslash \{ (i,I_j) \} \qquad \textrm{if} \ \alpha_i >1 \\
\I_< (\alpha , I ) \backslash \{ (i',I_{j'}) \mid i'=i \} \qquad \textrm{if} \ \alpha_i = 1. \\
\end{cases}$$
\end{proof}

\begin{definition}\label{def:sparseEN}
Adopt notation and hypotheses as in Setup \ref{set:linearStrandSetup}. Then the \emph{sparse Eagon-Northcott} complex (with respect to $<$) is the complex of Definition \ref{def:resnoflexideal}
\end{definition}

\begin{remark}
The complex of Definition \ref{def:resnoflexideal} is called \emph{sparse} because of the relationship of the differentials to that of the classical Eagon-Northcott complex. Namely, the differentials are obtained by simply setting some of the entries equal to $0$ in the matrix representation.
\end{remark}

\begin{cor}\label{cor:resnoflexinitialideal}
Adopt notation and hypotheses as in Setup \ref{set:linearStrandSetup}. Then the sparse Eagon-Northcott complex $E'_\bullet$ is a minimal free resolution of $R/\inn_< (I_n (M))$. 
\end{cor}

\begin{proof}
It is clear that the image of each basis element of $E_i'$ forms a linearly independent subset of $\ker (d_{i-1})$ for each $i \geq 2$. Using Theorem \ref{thm:boocherBnos}, this image must also be a generating set. 
\end{proof}

The following results are quick corollaries of Corollary \ref{cor:resnoflexinitialideal}. For the definition of the box polarization used in Corollary \ref{cor:boxPol}, see the introduction of \cite{polarizations}.

\begin{cor}
Adopt notation and hypotheses as in Setup \ref{set:linearStrandSetup}. Let $E_\bullet'$ denote the sparse Eagon-Northcott complex with respect to $<$. Then $E_\bullet' \otimes R / \sigma$ is a minimal free resolution of the quotient by the ideal of all squarefree monomials of degree $n$ in $m$ variables, where
\begin{align*}
\sigma = \{x_{11}-x_{21}, x_{11}-x_{31},\ldots,x_{11}-x_{n1}\} \cup \{x_{12}-x_{22},\ldots,x_{12}-x_{n2}\}\cup \ldots \\
\cup \{x_{1m}-x_{2m},\ldots, x_{1m}-x_{nm}\}.
\end{align*}
\end{cor}

\begin{proof}
This is immediate by Theorem \ref{thm:boocherBnos}.
\end{proof}

\begin{cor}\label{cor:boxPol}
Adopt notation and hypotheses as in Setup \ref{set:linearStrandSetup}. Let $E_\bullet'$ denote the sparse Eagon-Northcott complex with respect to $<$. Then, under the relabelling
$$x_{ij} \mapsto x_{j-i+1,i},$$
$E'_\bullet$ is a minimal free resolution of the quotient by the box polarization of $(x_1 , \dotsc , x_{m-n+1})^n$.

In particular, with the above relabelling, $E_\bullet' \otimes R/\sigma$ is a minimal free resolution of $(x_1 , \dotsc , x_{m-n+1})^n$, where
\begin{align*}
\sigma = \{x_{11}-x_{12}, x_{11}-x_{13},\ldots,x_{11}-x_{1n}\} \cup \{x_{21}-x_{22},\ldots,x_{21}-x_{2n}\}\cup \ldots \\
\cup \{x_{m-n+1,1}-x_{m-n+1,2},\ldots, x_{m-n+1,1}-x_{m-n+1,n}\}.
\end{align*}
\end{cor}

\section{Linear Strand of Certain DFIs}\label{sec:linStrandlcmclosed}
In this section, we construct an explicit \textit{linear strand} for the initial ideal of certain classes of DFIs. This linear strand is built as a \emph{generalized} sparse Eagon-Northcott complex, where the generators in each homological degree are parametrized by faces of the \emph{clique complex} associated to the simplicial complex $\Delta$ parametrizing the DFI. We then define the notion of an $i$-nonface of a simplicial complex (see Definition \ref{def:inonface}); it turns out that the nonexistence of $1$-nonfaces will guarantee that homology of the associated generalized sparse Eagon-Northcott complex is trivial in appropriate degrees. In particular, we deduce that so-called lcm-closed DFIs have the property that 
$$\beta_{i,i+n} (R/ J_\Delta) = \beta_{i,i+n} (R/ \inn_< J_\Delta) \ \textrm{for all} \ i \geq 0,$$
for any diagonal term order $<$. To conclude this section, we verify a conjecture of Ene, Herzog, and Hibi on the Betti numbers of binomial edge ideals when the associated graph has at most $2$ maximal cliques.

\begin{definition}\label{def:linStrand}
Let $F_\bullet$ be a minimal graded $R$-free complex with $F_0$ having initial degree $d$. Then the \emph{linear strand} of $F_\bullet$, denoted $F_\bullet^{\lin}$, is the complex obtained by restricting $d_i^F$ to $(F_i)_{d+i}$ for each $i \geq 1$.
\end{definition}

\begin{remark}
Observe that the minimality assumption in Definition \ref{def:linStrand} ensures that the linear strand is well defined. Choosing bases, the linear strand can be obtained by restricting to the columns where only linear entries occur in the matrix representation of each differential.
\end{remark}

Now we recall some results on linear strands of general $R$-modules due to Herzog, Kiani, and Madani.

\begin{theorem}[\cite{herzog2015linear}, Theorem 1.1]\label{thm:linstrandequiv}
Let $R$ be a standard graded polynomial ring over a field $k$. Let $G_\bullet$ be a finite linear complex of free $R$-modules with initial degree $n$. Then the following are equivalent:
\begin{enumerate}
    \item The complex $G_\bullet$ is the linear strand of a finitely generated $R$-module with initial degree $n$.
    \item The homology $H_i (G_\bullet)_{i+n+j} = 0$ for all $i>0$ and $j=0, \ 1$.
\end{enumerate}
\end{theorem}

\begin{prop}[\cite{herzog2015linear}, Corollary 1.2]\label{prop:linstrandcor}
Let $R$ be a standard graded polynomial ring over a field $k$. Let $G_\bullet$ be a finite linear complex of free $R$-modules with initial degree $n$ such that $H_i (G_\bullet)_{i+n+j} = 0$ for all $i>0$, $j=0, \ 1$. 

Let $N$ be a finitely generated $R$-module with minimal graded free resolution $F_\bullet$. Assume that there exist isomorphisms making the following diagram commute:
$$\xymatrix{G_1 \ar[d]_-{\sim} \ar[r] & G_0 \ar[d]^-{\sim} \\
F_1^{\textrm{lin}} \ar[r] & F_0^{\textrm{lin}}. \\}$$
Then $G_\bullet \cong F_\bullet^\textrm{lin}$.
\end{prop}

Notice that the following definition is simply a subcomplex of the sparse Eagon-Northcott complex as in Definition \ref{def:resnoflexideal}.

\begin{definition}\label{def:gensparseEN}
Adopt notation and hypotheses as in Setup \ref{set:linearStrandSetup}. Let $\Delta$ denote a simplicial complex on the vertex set $[m]$. Define $\cat{C}^<_0 (\Delta , M) := \bigwedge^n G$. For $i \geq 1$, let $\cat{C}^<_i (\Delta , M) \subseteq D_{i-1} (G^*) \otimes \bigwedge^{n+i-1} F$ denote the free submodule generated by all elements of the form
$$g^{* ( \alpha)} \otimes f_\sigma,$$
where $\sigma \in \Delta$ with $|\sigma|=n+i-1$ and $\alpha = (\alpha_1 , \dotsc , \alpha_n)$ with $|\alpha| = i-1$. 

Let $\cat{C}^<_\bullet (\Delta , M)$ denote the complex induced by the differential $$d'_1 : \bigwedge^n F \to \bigwedge^n G$$ sending $f_I \mapsto \inn_<  (M(f_I))$ and for $\ell\geq 2$, 
$$d'_\ell : D_{\ell - 1} (G^*) \otimes \bigwedge^{n+\ell-1} F \to D_{\ell-2} (G^*) \otimes \bigwedge^{n+\ell-2} F$$ is the sparse Eagon-Northcott differential
$$d_\ell ( g^{*(\alpha)} \otimes f_I ) = \sum_{\{i \mid \alpha_i >0\}} \sum_{(i,I_j) \in \I_< (\alpha , I)} (-1)^{j + 1}  x_{iI_j} g^{*(\alpha - \epsilon_i)} \otimes f_{I \backslash I_j }.$$
\end{definition}

\begin{definition}\label{def:genSparseENDef}
Adopt notation and hypotheses as in Setup \ref{set:linearStrandSetup}. Let $\Delta$ denote a simplicial complex on the vertex set $[m]$. The complex of Definition \ref{def:gensparseEN} will be called the \emph{generalized sparse Eagon-Northcott} complex.
\end{definition}
Notice that by the definition of a simplicial complex, the above complex is indeed well-defined. The following definition introduces a slight generalization of the notion of a \emph{minimal nonface}, that is, a set of vertices $\sigma\notin\Delta$ such that ever proper subset of $\sigma$ corresponds to a face of $\Delta$.

\begin{definition}\label{def:inonface}
Let $\Delta$ be a simplicial complex with vertex set $[m]$.
Then for some fixed $i$ and any $\ell\in \NN$, an $i$-nonface of $\Delta$ is given by an ordered list of vertices $\sigma = \{\sigma_1 < \dots < \sigma_\ell\}$ where each $\sigma_j \in [m]$
such that $\sigma \notin \Delta$ and for some $j\geq 1$ one has $\sigma \backslash \sigma_{j+k} \in \Delta$ for all $k = 0,\dots, i$. 
\end{definition}

\begin{remark} If $\Delta$ is an $(n-1)$-dimensional simplicial complex, we will often want to study $i$-nonfaces of cardinality at least $n+1$ of the clique complex $\clique$.
\end{remark}

\begin{example}\label{ex:exampleof1nonface}
Consider the following graph $G$ on vertices $\{1,2,3,4 \}$:
$$\xymatrix{3 \ar@{-}[d] & 2 \ar@{-}[dl] \\
4  & 1. \ar@{-}[l] \ar@{-}[ul] \ar@{-}[u] \\}$$
Observe that the associated clique complex has facets $\{1,2,4 \}$ and $\{1,3,4 \}$.
Then $\{1,2,3\}$ is a $1$-nonface of cardinality $3$ in $G^{\text{clique}}$, since both $\{1,2\}$ and $\{1,3\}$ are faces of $G^{\text{clique}}$; similarly, $\{2,3,4\}$ is also a $1$-nonface since both $\{3,4\}$ and $\{2,4\}$ are faces of the clique complex of $G$.
Similarly, $\{ 1,2,3,4 \}$ is a $1$-nonface of the clique complex of cardinality $4$, since $\{1,2,4 \}$ and $\{1,3,4 \}$ are both facets.

If we instead consider the graph
$$\xymatrix{3 \ar@{-}[d] & 2 \ar@{-}[l] \ar@{-}[dl] \\
4  & 1,  \ar@{-}[ul] \ar@{-}[u]  \\}$$
then the clique complex has facets $\{ 1,2,3\}$ and $\{2,3,4 \}$. 
In this case, there are no $1$-nonfaces of at least cardinality $3$ in the clique complex of $G$. For example, the set of vertices $\sigma = \{1,2,3,4\}$ is not a $1$-nonface because there are no two consecutive integers $i, i+1$ such that both $\sigma\setminus \sigma_i$ and $\sigma\setminus \sigma_{i+1}$ are faces of the associated clique complex.
Likewise, there are no $1$-nonfaces of cardinality $3$, so there are no $1$-nonfaces of at least cardinality $2$.

For the graph
$$\xymatrix{3 \ar@{-}[d] & 2 \ar@{-}[l] \ar@{-}[dl] \\
4  & 1, \ar@{-}[l] \ar@{-}[u]  \\}$$
the associated clique complex has facets $\{ 1,2,4\}$ and $\{2,3,4 \}$, and has no $1$-nonfaces of cardinality $4$.
However, $\{ 1,3,4 \}$ is a $1$-nonface of cardinality $3$ since $\{ 1,4 \}$ and $\{3,4 \}$ are vertices of $G$.
\end{example}

\begin{remark}
Notice that a minimal nonface $\sigma$ is a $\dim \sigma$-nonface in the above definition. Moreover, any $i$-nonface is a $k$-nonface for all $k \leq i$. 
\end{remark}

In the proofs of the results in the remainder of this section, notice that we have chosen to augment our complexes with the ring $R$. This means that we are resolving the quotient ring $R/I$ as opposed to the module $I$; this has the effect of shifting the indexing in the statements of Theorem \ref{thm:linstrandequiv} and Proposition \ref{prop:linstrandcor}.

\begin{lemma}\label{lem:linstrandnonfaces}
Adopt notation and hypotheses as in Setup \ref{set:linearStrandSetup}. If the simplicial complex $\Delta$ has no $1$-nonfaces of cardinality $\geq n+1$, then the complex $\cat{C}^<_\bullet (\Delta , M)$ is the linear strand of a finitely generated graded $R$-module with initial degree $n$.
\end{lemma}

\begin{proof}
We employ Theorem \ref{thm:linstrandequiv}. To avoid trivialities, assume $n \leq \dim (\Delta) + 1$. Observe first that $H_i (\cat{C}_\bullet^< ( \Delta , M))_{i+n-1} = 0$ for all $i \geq 1$ trivially.

To finish the proof, it will be shown that
$$H_i (\cat{C}_\bullet^< ( \Delta , M))_{i+n} \neq 0 \ \textrm{for all} \ i \geq 1$$
$$\implies$$
$$\textrm{there exists a} \ 1\textrm{-nonface of cardinality} \ n+i, \ \textrm{for all} \ i\geq 1.$$
For convenience, use the notation $\cat{C}_\bullet^< (M) =: E_\bullet'$, where $\cat{C}_\bullet^< (M)$ is as in Definition \ref{def:gensparseEN}. 

Assume $H_i (\cat{C}_\bullet^< ( \Delta , M))_{i+n} \neq 0$. Let $z \in \cat{C}_i^< (\Delta , M)$ be a cycle that is not a boundary; without loss of generality, assume $z$ is multihomogeneous. The complex $\cat{C}_\bullet^< (M)$ is exact, whence $z = d(y)$ for some $y \in \cat{C}_{i+1}^< (M)$. By multihomogeneity, $y = \lambda g^{*(\alpha)} \otimes f_\sigma$ for some $\lambda \in k^\times$ with $|\sigma|=n+i$, $|\alpha|=i$. The assumption that $z$ is not a boundary implies that $\sigma \notin \Delta$, since otherwise $y \in \cat{C}_{i+1} (\Delta, M)$. By definition of the differential of $\cat{C}_\bullet^< (M)$,
$$z = \lambda \cdot \sum_{\{\ell \mid \alpha_\ell >0\}} \sum_{(\ell,\sigma_j) \in \I_< (\alpha , \sigma)} (-1)^{j + 1} x_{\ell \sigma_j} g^{*(\alpha - \epsilon_\ell)} \otimes f_{\sigma \backslash \sigma_j }.$$
Since $z \neq 0$, one has $(\ell , \sigma_j) \in \I_< (\alpha , \sigma)$ for some $\ell$, $j$. By definition of $\I_< (\alpha , \sigma)$, this means $(\ell , \sigma_k) \in \I_< ( \alpha , \sigma)$ for all $|\alpha_{\leq \ell-1} | \leq j \leq |\alpha_{\leq \ell}|$. This translates to the fact that $\sigma$ is an $\alpha_\ell$-nonface of cardinality $n+i$. Since $\alpha_\ell \geq 1$, the result follows.
\end{proof}

\begin{remark}
The proof of Lemma \ref{lem:linstrandnonfaces} allows one to construct explicit examples of nonzero homology on the complex $\cat{C}_\bullet^< (\Delta , M)$. Let $\clique$ be the simplicial complex associated to the first graph of Example \ref{ex:exampleof1nonface}. Then the element
$$z = x_{22} f_{1,3,4} - x_{23} f_{1,2,4}$$
is a cycle which is not a boundary in $\cat{C}_\bullet^< ( \clique ,M)$. 
\end{remark}

\begin{lemma}\label{lem:1sthomologylemma}
Adopt notation and hypotheses as in Setup \ref{set:linearStrandSetup}. Then the following are equivalent:
\begin{enumerate}
    \item $H_1 ( \cat{C}_\bullet^< (\Delta , M))_{n+1} = 0$,
    \item $\Delta$ has no $1$-nonfaces of cardinality $n+1$.
\end{enumerate}
\end{lemma}

\begin{proof}
The implication $(2) \implies (1)$ is Lemma \ref{lem:linstrandnonfaces}. Conversely, assume that $\sigma \notin \Delta$ is a $1$-nonface of cardinality $n+1$. By definition, there exists some $j$ such that $\sigma \backslash \sigma{j}$ and $\sigma \backslash \sigma_{j+1} \in \Delta$. This means that  $z= (-1)^{j+1} (x_{j\sigma_j} f_{\sigma \backslash \sigma_j} -  x_{j\sigma_{j+1}} f_{\sigma \backslash \sigma_{j+1}})$ is a cycle in $\cat{C}_1^< (\Delta ,M)$ that is not a boundary, since $z = d_2 (g_j^* \otimes f_\sigma )$, and $g_j^* \otimes f_\sigma \notin \cat{C}_2 (\Delta , M)$ by construction.  
\end{proof}

Recall that the standard Eagon-Northcott complex inherits a $\ZZ^n \times \ZZ^m$-grading, as described in Notation \ref{not:ENgrading}. Since the sparse Eagon-Northcott complexes of Definition \ref{def:gensparseEN} are obtained by simply setting certain entries in the differentials equal to $0$, these maps will remain multigraded in an identical manner. We tacitly use this multigrading for the remainder of this section. 

\begin{theorem}\label{thm:linearstrandforclosedDelta}
Adopt notation and hypotheses as in Setup \ref{set:linearStrandSetup}. Assume that $\Delta$ is an $(n-1)$-pure simplicial complex such that $\clique$ has no $1$-nonfaces of cardinality $n+1$. Let $F_\bullet$ denote the minimal graded free resolution of $\inn_<(J_\Delta)$; then
$$F_\bullet^\textrm{lin} \cong \cat{C}_\bullet^< (\clique, M).$$
\end{theorem}

\begin{proof}
Let $Z^\textrm{lin} := (\ker d_1)_{n+1}$, where $d_1$ is the first differential of the complex $\cat{C}_\bullet^< (\clique , M)$. By construction, $\cat{C}_1^< (\clique , M)$ is generated in degree $n+1$ and hence induces a homogeneous map
$$\partial : \cat{C}_1^< ( \clique , M) \to Z^\textrm{lin}.$$
Let $0 \neq z \in Z^{\textrm{lin}}$ be an element of multidegree $(\alpha_1, \alpha_2) = (\epsilon_s + \one,\epsilon_{i_1}+\cdots + \epsilon_{i_{n+1}})\in \ZZ^n\times \ZZ^m$ (where $\one$ denotes the appropriately sized vector of all $1$'s). 
Set $\tau := \{ i_1 < \cdots < i_{n+1} \}$ corresponding to ordering the elements in the support of $\alpha_2$; by multihomogeneity, there are constants $\lambda_k \in k$ such that
$$z = \sum_{k=1}^{n+1} \lambda_k x_{s i_k} f_{\tau \backslash i_k}.$$
Since $z$ is a cycle of $\cat{C}_1^< (M)$ (where $\cat{C}_\bullet^< (M) := E_\bullet'$ is as in Definition \ref{def:gensparseEN}), there exists $y \in \cat{C}_2^< (M)$ such that $d_2(y) = z$. By multihomogeneity, $y = \lambda g_s \otimes f_\tau$ for some constant $\lambda$, whence $z = \lambda (-1)^{s+1} (x_{s\sigma_s} f_{\sigma \backslash \sigma_s} - x_{s\sigma_{s+1}} f_{\sigma \backslash \sigma_{s+1}} )$. This implies that $\sigma \in \clique$, since otherwise $\clique$ would have a $1$-nonface of cardinality $n+1$, contradicting our assumptions. Thus $Z^\textrm{lin}$ is generated by the set
$$\{r_s (\sigma) := (-1)^{s+1} (x_{s\sigma_s} f_{\sigma \backslash \sigma_s} - x_{s\sigma_{s+1}} f_{\sigma \backslash \sigma_{s+1}} ) \mid 1 \leq s \leq n, \ \sigma \in \clique, \ |\sigma| = n+1 \}.$$
Moreover, since $\textrm{mdeg} (r_s(\sigma)) \neq \textrm{mdeg} (r_{s'} (\sigma')$ for $s\neq s'$ or $\sigma \neq \sigma'$, the above is a basis. Finally, $d_2 (g_s^* \otimes f_\sigma ) = r_s (\sigma)$, whence the induced map $\partial$ is an isomorphism of vector spaces.
\end{proof}

\begin{remark}
Let $\Delta$ be an $(n-1)$-pure closed simplicial complex. Then $\clique$ has no minimal nonfaces in cardinality $\geq n+1$, since any minimal nonface is in particular a $1$-nonface. This means that $\clique$ satisfies the hypotheses of Theorem $3.1$ of \cite{herzog2015linear}. 
\end{remark}

The following theorem says that the Betti numbers in the linear strand of a DFI and its initial ideal coincide in the case that the associated clique complex has no $1$-nonfaces. Recall that in general, the Betti numbers of the initial ideal are merely an upper bound and it is quite rare to have equality throughout.

\begin{theorem}\label{thm:mainpunchline}
Adopt notation and hypotheses as in Setup \ref{set:linearStrandSetup}. Assume that $\Delta$ is an $(n-1)$-pure simplicial complex with no $1$-nonfaces of cardinality $n+1$. Then for all $i \geq 1$,
$$\beta_{i,n+i} (J_\Delta) = \beta_{i,n+i} (\inn_< (J_\Delta)).$$
\end{theorem}

\begin{proof}
Notice that the linear strand of $J_\Delta$ is $\cat{C}_\bullet (\clique , M)$ where $\cat{C}_\bullet$ is the generalized Eagon-Northcott complex of \cite{herzog2015linear}. Then, $\cat{C}_\bullet$ and $\cat{C}_\bullet^<$ have the same underlying free modules, so the result follows.
\end{proof}

\begin{prop}\label{prop: no1nonface} Let $\Delta$ be pure $(n-1)$-dimensional simplicial complex on $m$ vertices, and let $J_\Delta$ be its associated $n$-determinantal facet ideal. If $J_\Delta$ is lcm-closed, then there are no $1$-nonfaces in $\clique$.
\end{prop}

\begin{proof}
It suffices to show that there are no $1$-nonfaces of cardinality $n+1$ in $\clique$. Suppose, seeking contradiction, that $\ff = \{f_1 < \dots < f_{n+1}\}$ is a $1$-nonface of cardinality $n+1$ in $\clique$. By definition, there exists some $f_i$ such that $\aa = \{f_1,\dots, \widehat{f_i}, f_{i+1}, \dots, f_{n+1}\}$ and $\bb = \{f_1,\dots, f_i, \widehat{f_{i+1}},\dots, f_{n+1}\}$ are facets of $\Delta$. Let $\Delta_\aa$ and $\Delta_\bb$ be maximal cliques of $\Delta$ containing $\aa$ and $\bb$, respectively, with nontrivial intersection.
By assumption, $\Delta_\aa\neq \Delta_\bb$. Because $J_\Delta$ is lcm-closed, there exists some facet $\cc\in \Delta_\aa\cap \Delta_\bb$ such that $\inn[\cc]$ divides
$$\lcm (\inn[\aa],\inn[\bb]) = x_{1 f_1} x_{2 f_2} \dots x_{i f_i} x_{i f_{i+1}} x_{i+1 f_{i+2}} \dots x_{n,f_{n+1}}.$$
The only possible facets $\cc$ of $\Delta$ satisfying this property are $\aa$ and $\bb$ themselves, so they must be in the same maximal clique of $\Delta$ and $\ff$ is indeed a face in $\clique$, giving the desired contradiction.
\end{proof}



\begin{cor}\label{cor:lcmClosed}
Adopt notation and hypotheses as in Setup \ref{set:linearStrandSetup}. Assume that $J_\Delta$ is an lcm-closed DFI. Then for all $i$,
$$\beta_{i,n+i} (J_\Delta) = \beta_{i,n+i} (\inn_< (J_\Delta)).$$
\end{cor}

Using Corollary \ref{cor:lcmClosed} as initial evidence, we pose the following:

\begin{conj}\label{conj:BnosConj}
Adopt notation and hypotheses as in Setup \ref{set:linearStrandSetup}. Assume that $J_\Delta$ is an lcm-closed DFI. Then
$$\beta_{ij} (R/J_\Delta) = \beta_{ij} (R / \inn_< J_\Delta) \quad \textrm{for all} \ i,j.$$
\end{conj}

\begin{remark}
It is important to notice that the condition on $1$-nonfaces is not sufficient for the minimal generators to form a Gr\"obner basis, and is hence more general than the property of being lcm-closed. For example, let $\Delta$ be the simplicial complex with facets ${1,2,3}$, ${1,4,5}$, and ${1,6,7}$. One can verify directly that there are no $1$-nonfaces of cardinality $4$, but calculations in Macaulay2 show that the minimal generating set of the associated determinantal facet ideal does \emph{not} form a Gr\"obner basis.
\end{remark}

To conclude this section, we gather some necessary facts to prove (in a special case) a conjecture of Ene, Herzog, and Hibi. The following result, due to Conca and Varbaro, shows that ideals with squarefree initial ideals (with respect to some term order) exhibit homological behavior similar to that of the associated generic initial with respect to revlex. Recall that a Betti number $\beta_{i,j}(N)\neq 0$ of an $R$-module $N$ is \textit{extremal} if $\beta_{l,r} = 0$ for all $l\geq i$, $r\geq j+1$, and $r-l\geq j-i$; that is, $\beta_{ij}$ is the nonzero top left ``corner'' in a block of zeroes in the Betti diagram of $N$.

\begin{theorem}[{\cite[Follows from Theorem 1.3]{conca2020square}}]\label{thm:concavarbaro}
Let $I$ be a homogeneous ideal in a standard graded polynomial ring $R$ with term order $<$. If $\inn_< (I)$ is squarefree, then the extremal Betti numbers of $R/I$ and $R/ \inn_< (I)$ coincide. In particular,
$$\reg (R / I) = \reg (R / \inn_< (I)) \quad \textrm{and} \quad \pd_R (R/I) = \pd_R (R/ \inn_< (I)).$$
\end{theorem}

We recall the definition of a \textit{consecutive cancellation}, introduced by Peeva in \cite{peeva2004consecutive}.
\begin{definition}
A sequence $\{q_{i,j}\}$ of numbers is obtained from a sequence $\{p_{i,j}\}$ by a \textit{consecutive cancellation} if there exist indices $s$ and $r$ such that
$$
q_{s,r} = p_{s,r}-1, \quad q_{s+1, r} = p_{s+1, r} - 1,
$$
and $q_{i,j} = p_{i,j}$ for all other values of $i$ and $j$.
\end{definition}

\begin{theorem}[{\cite[Theorem 1.1]{peeva2004consecutive}}]\label{thm:consecutiveCancellation}
Let $I$ be a graded ideal and let $<$ be any term order. Then the graded Betti numbers $\beta_{i,j} (R/I)$ can be obtained from the graded Betti numbers $\beta_{i,j} (R/ \inn_< (I))$ by a sequence of consecutive cancellations.
\end{theorem}

\begin{theorem}[{\cite[Corollary of Theorem 2.3]{malayeri2020proof}}]\label{thm:regularitybound} 
Let $G$ be a closed graph with at most $2$ maximal cliques. Then
$$\reg (R / J_G) \leq 2.$$
\end{theorem}
Finally, we arrive at the last result of this section.

\begin{cor}\label{cor:EHHconj}
Let $G$ be a closed graph with at most $2$ maximal cliques. Let $J_G$ denote the associated binomial edge ideal $J_G$ and let $<$ denote any diagonal term order. Then for all $i, j$,
$$\beta_{i,j} (R / J_G) = \beta_{i,j} ( R / \inn_< ( J_G)).$$
\end{cor}

\begin{proof}
Since $G$ is a closed graph, $J_G$ is lcm-closed and hence the associated clique complex of $G$ has no $1$-nonfaces (by Proposition \ref{prop: no1nonface}). It is well known that every binomial edge ideal has squarefree Gr\"obner basis with respect to $<$ (see \cite{herzog2010binomial}); in particular, Theorem \ref{thm:concavarbaro} conbimed with Theorem \ref{thm:regularitybound} shows that $\reg (R / \inn_< J_G) \leq 2$. Theorem \ref{thm:consecutiveCancellation} asserts that the Betti numbers of $R/ J_G$ are obtained by those of $R / \inn_< J_G$ by a sequence of consecutive cancellations. However, Theorem \ref{thm:mainpunchline} implies that no cancellations are possible.
\end{proof}

%% file: sections/cellular.tex
\section{Linear Strands Supported on Cellular Complexes}\label{sec: cellular}

In this section, we introduce the notion of a linear strand supported on a polyhedral cell complex (Definition \ref{def: cellLinStrand}), generalizing the well-studied phenomenon of cellular resolutions.
We show that the first linear strand of the initial ideal of any determinantal facet ideal with respect to a diagonal term order is supported on a polyhedral cell complex. In particular, this cell complex is an induced subcomplex of the \emph{complex of boxes} introduced by Nagel and Reiner (see Construction \ref{const: complexOfBoxes}), which supports a minimal linear free resolution of squarefree strongly stable and strongly stable ideals \cite{NR09}.

We begin by recalling some basic notions from the theory of cellular resolutions. For a more detailed exposition, see \cite[Chapter 4]{miller2004combinatorial}.

\begin{definition}\label{def: cellComplex} A \emph{polyhedral cell complex} $\cP$ is a finite collection of convex polytopes (called \emph{cells} or \emph{faces} of $\cP$) in some Euclidean space, satisfying the following two properties:
\begin{itemize}
    \item if $H$ is a polytope in $\cP$, then every face of $H$ also lies in $\cP$, and
    \item if $H_i, H_j$ are both in $\cP$, then $H_i\cap H_j$ is a face of both $H_i$ and $H_j$.
\end{itemize}
Denote by $V(\cP)$ the set of vertices (or $0$-dimensional cells) of $\cP$.
If $\cX\subseteq V(\cP)$, the \emph{induced subcomplex} of $\cP$ on $\cX$ is the subcomplex $\{F\in \cP \mid V(F)\subseteq \cX \}$.
The $f$-vector of a $d$-dimensional polyhedral cell complex $\cP$ is the vector $(f_0,f_1,\dots, f_{d} )$, where $f_i$ is the number of $i$-dimensional cells of $\cP$. Any polyhedral cell complex is equipped with a reduced cellular chain complex $(C_\doot, \partial_\doot)$, where the signs in the differentials are specified by (arbitrarily) orienting the faces of $\cP$.
\end{definition}

\begin{construction}\label{const: cellRes} Set $S = k[x_1,\dots, x_n]$ to be a polynomial ring over a field $k$.
Let $\cP$ be an oriented polyhedral complex and let $(\alpha_H)_{H\in\cP}\subseteq \ZZ^{n}$ be a labeling of the cells of $\cP$ such that
$$
\xx^{\alpha_H} = \lcm\{\xx^{\alpha_G} \mid G\subset H\}.
$$
The labeled complex $(\cP,\alpha)$ gives rise to an algebraic complex of free $\ZZ^n$-graded $S$-modules in the following way. 
For two cells $G,H\in \cP$ with $\dim H = \dim G + 1$ denote by $\epsilon(H,G)\in \{0,\pm 1\}$ the coefficient of $G$ in the cellular boundary of $H$. Define the free modules
$$
F_i \coloneqq \bigoplus_{\substack{H\in \cP\\ \dim H = i+1}} S(-\alpha_H).
$$
The differentials $d_i: F_i\ra F_{i-1}$ are given by
$$
d(e_H) \coloneqq \sum_{\dim G = \dim H -1} \epsilon(H,G) \xx^{\alpha_H-\alpha_G} e_G.
$$
One can verify this defines an algebraic complex $\cF_\cP$. For $\beta\in \ZZ^n$, denote by $\cP_{\leq \beta}$ the subcomplex of $\cP$ consisting of all cells $H\in \cP$ with $\alpha_H\leq \beta$ coordinatewise.
Let $I_{\cP} = \langle \xx^{\alpha_v} \mid v\in \cP \text{ a vertex}\rangle$.
\end{construction}

\begin{lemma}[{\cite[Proposition 4.5]{miller2004combinatorial}}]\label{lem: cellRes} Adopt notation and hypotheses of Construction \ref{const: cellRes}. Let $\cF_\cP$ be the algebraic complex obtained from the labeled polyhedral complex $(\cP,\alpha)$. If for every $\beta\in \ZZ^n$ the subcomplex $\cP_{\leq \beta}$ is acyclic over $k$, then $\cF_\cP$ resolves the quotient $S/{I_{\cP}}$. Furthermore, the resolution is minimal if $\alpha_H\neq \alpha_G$ for any two faces $G\subset H$ with $\dim H = \dim G+1$.
\end{lemma}

\begin{definition}\label{def: cellRes} Adopt notation and hypotheses of Construction \ref{const: cellRes}. The complex $\cF_\cP$ is a \emph{cellular resolution} if it meets the criteria of Lemma \ref{lem: cellRes}, and the polyhedral complex $\cP$ \emph{supports} the resolution.
\end{definition}

We extend the notion of cellular resolution to study multigraded linear strands supported on a polyhedral cell complex. The following lemma follows naturally from Theorem \ref{thm:linstrandequiv}.

\begin{lemma}\label{lem: cellLinStrand} Assume the complex $\cF_\cP$ is $d$-linear, i.e., all the generators of $I_{\cP}$ have degree $d$ and all higher syzygy maps are given by linear forms. Then $\cF_\cP$ is the first linear strand of $S/{I_{\cP}}$ if, for any $\alpha\in \ZZ^n$ with $\abs{\alpha} = d+k$ and $k>0$, $\tilde H_i(\cP_{\leq \alpha}) = 0$ for $i = k$ and $k-1$.
\end{lemma}

\begin{definition}\label{def: cellLinStrand} The complex $\cF_\cP$ is a \emph{cellular linear strand} if it satisfies the criteria of Lemma \ref{lem: cellLinStrand}. In this case, we say that the polyhedral complex $\cP$ \emph{supports the linear strand} of $S/{I_\cP}$.
\end{definition}


The following construction by Nagel and Reiner can be defined more generally for squarefree strongly stable and strongly stable ideals, but for our purposes we only consider the case when the ideal in question is $\inn(I_n(M))$.

\begin{construction}\label{const: complexOfBoxes}(see \cite{NR09})
Partition the variables of $R$ into $n$ subsets $\Xv_i = \{x_{i1},\dots, x_{im}\}$. Set $K = \{\aa \mid \aa \text{ is an $n$ subset of } [m]\}$, so the elements of $K$ are in one-to-one correspondence with the generators of $\inn_<(I_n(M))$ via $\aa = \{a_1 < \dots < a_n\} \longleftrightarrow \xx_\aa = x_{1 a_1} \cdots x_{n a_n}$.

Call a subset of $K$ which is a Cartesian product $X_1 \times \dots \times X_n$ for subsets $X_j\subseteq \Xv_j$ a \emph{box} inside $K$, and define the \emph{complex of boxes} of $K$ to be the polyhedral subcomplex of the product of simplices $2^{\Xv_1}\times \dots \times 2^{\Xv_n}$ having faces indexed by the boxes inside $K$.
\end{construction}

\begin{theorem}[{\cite[Theorem 3.12]{NR09}}]\label{thm: cplxBoxesRes}
Labelling a vertex $\aa$ in the complex of boxes by the monomial $\xx_\aa$ gives a minimal linear cellular resolution of $R/\inn_<(I_n(M))$, where $<$ is any diagonal term order.

Hence $\beta_{i,j}(R/\inn(I_n(M))) = 1$ where $i = \sum_k \abs{X_k}-n$, $j = X_1\sqcup\dots\sqcup X_n$ for every box $X_1\times \dots \times X_n$ inside $K$, and all other Betti numbers vanish.
\end{theorem}

\begin{notation}\label{not: subCplx}
Let $\Delta$ be a pure $(n-1)$-dimensional simplicial complex on $m$ vertices, and let $\cP$ denote the complex of boxes supporting a minimal linear cellular resolution of $R/\inn(I_n(M))$. Denote by $\cP(\Delta)$ the induced polyhedral subcomplex of $\cP$ on the vertex set labeled by $\{\xx_\aa \mid \aa \text{ a facet of } \Delta \}.$
\end{notation}

\begin{theorem}\label{thm: cplxBoxesLinStrand} Let $\Delta$ be a pure $(n-1)$-dimensional simplicial complex on $m$ vertices. Then the linear strand of $R/\inn(J_\Delta)$ is supported on $\cP(\Delta)$.
\end{theorem}

\begin{proof}
First, observe that although $\inn(J_\Delta)$ may have generators in higher degree, its linear strand is completely determined by syzygies among the generators of degree $n$.

Let $\ell = f_{n-1}(\Delta)$, the number of facets of $\Delta$, and proceed by downward induction on $\ell$.
The base case $\ell = \binom{m}{n}$ corresponds to the case where $J_\Delta = I_n(M)$ and is clear. Fix $\ell\geq 0$ and assume that for any $\Delta$ with $\ell$ generators, $\cP(\Delta)$ supports the first linear strand of $\inn(J_\Delta)$. Let $\Delta'\subset \Delta$ be the subcomplex with a single facet $\aa$ removed. Then $\cP(\Delta')\subset \cP(\Delta)$ is the induced subcomplex $\cP(\Delta)\setminus v_\aa$, where $v_\aa$ is the vertex labeled by the generator $\xx_\aa$ of $\inn(J_\Delta)$.

By Lemma \ref{lem: cellLinStrand}, it suffices to check that for any multidegree $\alpha$ with $\abs{\alpha} = n+k$, $\tilde H_k(\cP(\Delta')_{\leq \alpha}) = \tilde H_{k-1}(\cP(\Delta')_{\leq \alpha}) = 0$.
Observe that any face of dimension $k$ in $\cP(\Delta')_{\leq \alpha}$ has multidegree $\beta$ such that $\abs{\beta} = \abs{\alpha} = n+k$. Therefore, if $\dim \cP(\Delta')_{\leq \alpha} = k$, it is the unique box in $\cP$ with multidegree $\alpha$ by Theorem \ref{thm: cplxBoxesRes}, which is contractible. If $\dim \cP(\Delta')_{\leq \alpha}<k$, then $H_k(\cP(\Delta')_{\leq \alpha}) =0 $ trivially.

By the inductive hypothesis, $\tilde H_{k-1}(\cP(\Delta')_{\leq \alpha}) = 0$ as long as $\xx_\aa$ does not divide $\xx^\alpha$, so suppose it does. Let $C$ be a cycle of dimension $k-1$ in $\cP(\Delta')_{\leq \alpha}$.
Since $\tilde H_{k-1}(\cP(\Delta)_{\leq \alpha}) = 0$ by the inductive hypothesis, there is some boundary $B$ of dimension $k$ in $\cP(\Delta)_{\leq \alpha}$ of degree $\alpha$ containing the vertex $v_\aa$. By Theorem \ref{thm: cplxBoxesRes}, there is a unique box in $\cP$ with multidegree $\alpha$, so this must be $B$. But then $\del(B)$ will be a linear combination of its codimension $1$ faces, including those containing $v_\aa$, so it cannot be $C$.
\end{proof}

As a consequence of Theorem \ref{thm: cplxBoxesLinStrand}, we obtain the following means of computing the Betti numbers in the first linear strand of $\inn(J_\Delta)$.

\begin{cor}\label{cor: fVectorBettis} Let $\Delta$ be a pure $(n-1)$-dimensional simplicial complex on $m$ vertices and let $\cP(\Delta)$ be as in Notation \ref{not: subCplx}. Then $\beta_{i,i+1}(R/\inn(J_\Delta)) = f_{i}(\cP(\Delta))$.
\end{cor}

\input{figures/boxCplx}
\input{figures/linStrand}

\begin{example}
Let $G$ be the graph in Example \ref{ex:exampleof1nonface} with clique decomposition $\{1,2,4\} \cup \{1,3,4\}$. The complex of boxes $\cP$ for $\inn(I_2(M))$ where $M$ is a $2\times 4$ matrix is shown in Figure \ref{fig: boxCplx}.
The induced subcomplex of $\cP(G)$ of $\cP$ on the vertices corresponding to edges in $G$ is depicted in Figure \ref{fig: linStrand}. The $f$-vector for this subcomplex is $(5,6,2)$, which indeed corresponds to the Betti numbers in the linear strand of $\inn(J_G)$. In particular, $\cP(G)$ has $1$-cells corresponding to the $1$-nonfaces $\{1,2,3\}$ and $\{2,3,4\}$ as well as $2$-cells corresponding to the $1$-nonface $\{1,2,3,4\}$ in Example \ref{ex:exampleof1nonface}.

Note, however, that this cell complex cannot support the full free resolution of the ideal generated by the labels on its vertices, since $\cP(\Delta)_{\leq \alpha}$ is not in general acyclic for any multidegree $\alpha$. Consider, for example, the multidegree $\alpha = x_{11} x_{12} x_{23} x_{24}$. Then $\cP(G)_{\leq \alpha}$ consists of the disjoint vertices labeled by $x_{11}x_{23}$ and $x_{12}x_{24}$, so $\tilde H_0(\cP(G)_{\leq \alpha}) = 1$ and this complex does not satisfy the hypotheses of Lemma \ref{lem: cellRes}. 
\end{example}

\begin{remark}
Nagel and Reiner show that other strongly stable and squarefree strongly stable ideals have a minimal, linear, cellular resolution given by the complex of boxes $\cP$. However, ideals generated by a subset of generators in these other cases do not, in general, have a linear strand supported on the induced subcomplex of $\cP$ in the same manner. The key issue is that in other cases, the multigraded Betti numbers may not be $0$ or $1$, so the proof of Theorem \ref{thm: cplxBoxesLinStrand} does not apply.
\end{remark}

%% file: figures/boxCplx.tex
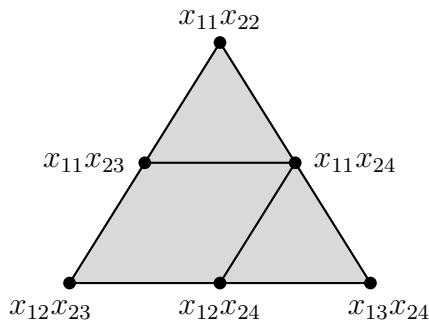
\begin{figure}[h]
\begin{tikzpicture}[]
\draw[fill=gray!30]   (-2,0) -- (2,0) -- (0,3.2) -- (-2,0);
\filldraw[thick]
(-2,0) circle (2pt) --
(0,0) circle (2pt)  -- 
(2,0) circle (2pt) --
(1,1.6) circle (2pt) --
(0,3.2) circle (2pt) --
(-1,1.6) circle (2pt) --
(-2,0);
\draw[thick] (-1,1.6)--(1,1.6)--(0,0);
\node at (-2.25,-0.35) {$x_{12} x_{23}$};
\node at (0,-0.35) {$x_{12} x_{24}$};
\node at (2.25,-0.35) {$x_{13} x_{24}$};
\node at (-1.8,1.6) {$x_{11} x_{23}$};
\node at (1.8, 1.6) {$x_{11} x_{24}$};
\node at (0, 3.5) {$x_{11} x_{22}$};
\end{tikzpicture}
\caption{Complex of boxes $\cP$ for $\inn(I_2(M))$ where $M$ is a $2\times 4$ matrix.}\label{fig: boxCplx}
\end{figure}

%% file: figures/linStrand.tex
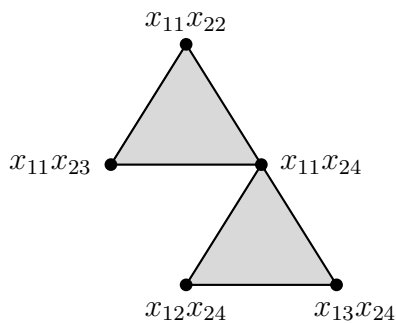
\begin{figure}[h]
\begin{tikzpicture}[]
\draw[fill=gray!30] (1,1.6)--(-1,1.6)--(0,3.2);
\draw[fill=gray!30] (1,1.6)--(0,0)--(2,0);
\filldraw[thick]
(0,0) circle (2pt)  -- 
(2,0) circle (2pt) --
(1,1.6) circle (2pt) --
(0,3.2) circle (2pt) --
(-1,1.6) circle (2pt);
\draw[thick] (-1,1.6)--(1,1.6)--(0,0);
\node at (0,-0.35) {$x_{12} x_{24}$};
\node at (2.25,-0.35) {$x_{13} x_{24}$};
\node at (-1.8,1.6) {$x_{11} x_{23}$};
\node at (1.8, 1.6) {$x_{11} x_{24}$};
\node at (0, 3.5) {$x_{11} x_{22}$};
\end{tikzpicture}
\caption{$\cP(G)$ where $G$ has clique decomposition $\{1,2,4\}\cup \{1,3,4\}$ as in Example \ref{ex:exampleof1nonface}.}\label{fig: linStrand}
\end{figure}